\documentclass[12pt]{article}
\usepackage[amsmath]{e-jc}
\usepackage{tikz}
\usetikzlibrary{arrows.meta, decorations.markings}
\usepackage{latexsym}
\usepackage{mathrsfs}
\usepackage{wasysym}
\usepackage{float}
\usepackage{graphicx}
\usepackage{calligra}
\allowdisplaybreaks[4]

\dateline{May 12, 2024}{July 28, 2025}{TBD}


\Copyright{The authors. Released under the CC BY-ND license (International 4.0).}

\title{Paths with two blocks in oriented graphs\\ of large minimum semi-degree}

\author{Bin Chen\authornote{1}\and Xinmin Hou\authornote{2,3}\and Xinyu Zhou\authornote{2}}

\authortext{1}{School of Mathematics and Statistics, Fuzhou University, Fujian, China.}
\authortext{2}{School of Mathematical Sciences, University of Science and Technology of China, Hefei 230026, Anhui, China.}
\authortext{3}{CAS Key Laboratory of Wu Wen-Tsun Mathematics, University of Science and Technology of China, Hefei 230026, Anhui, China.}

\begin{document}
\maketitle


\begin{abstract}
	
	Stein (2020) conjectured that for any positive integer $k$, every oriented graph of minimum semi-degree greater than $k/2$ contains every oriented path of length $k$.
	This conjecture is true for directed paths by a result from Jackson (JGT, 1981).
	In this paper, we establish the validity of Stein's conjecture specifically for any oriented path with two blocks, where, a block of an oriented path $P$ refers to  a maximal directed subpath within $P$.
	
\end{abstract}

\section{Introduction}

\noindent

Throughout this paper, all digraphs considered have no loops or parallel arcs. Let $D=(V, A)$ be a digraph with vertex set $V=V(D)$ and arc set $A=A(D)$. For an arc  $(u, v)\in A(D)$, we call $u$ the \emph{tail} and $v$ the \emph{head} of the arc $(u,v)$, and we also say $u$ {\em dominates} $v$ or $v$ is {\em dominated} by $u$. Let $\delta^{+}(D)$ and $\delta^{-}(D)$ be the \emph{minimum out-degree} and \emph{minimum in-degree} of $D$, respectively. By $\delta^{0}(D)$, we denote the \emph{minimum semi-degree} of $D$, which is the smallest value between $\delta^{+}(D)$ and $\delta^{-}(D)$. Two different vertices are \emph{adjacent} if there is an arc connecting them, and two different arcs are adjacent if they share a common vertex.

An \emph{orientation} of a graph $G$ is an assignment of exactly one direction to every edge of $G$. A \emph{directed path}\;(resp., \emph{anti-directed path}) is an oriented path such that any two adjacent arcs have the same orientation\;(resp., opposite orientations). A \emph{directed cycle}\;(resp., \emph{anti-directed cycle}) is an oriented cycle such that every pair of two adjacent arcs have the same orientation (resp., opposite orientations). A \emph{Hamiltonian oriented path}\;(resp., \emph{Hamiltonian oriented cycle}) is an oriented path\;(resp., oriented cycle) that passes through each vertex of a given digraph.

A digraph is \emph{traceable}\;(resp., \emph{anti-traceable}) if there is a Hamiltonian directed path\;(resp., Hamiltonian anti-directed path).
A digraph is \emph{Hamiltonian} if it contains a Hamiltonian directed cycle.
A digraph $D$ is \emph{strongly connected}\;(\emph{strong} for short) if and only if for every pair of distinct vertices $u, v\in V(D)$, there is a directed path from $u$ to $v$ and a directed path from $v$ to $u$.

Paths and cycles are one of the most fundamental objects in graph theory. Given a directed graph, a natural question is to find sufficient conditions for the existence of oriented paths or  cycles of specific types and lengths  within this digraph.

R\'{e}dei \cite{R} proved that any tournament is traceable and contains an odd number of Hamiltonian directed paths. In \cite{C}, Camion showed that every strong tournament is Hamiltonian. Extending Camion's result, Moon \cite{M} further proved that every strong tournament is vertex pancyclic. Let $RT_{n}$ be a regular tournament on $n$ vertices and $PT_{n}$ be a Paley tournament on $n$ vertices (A Paley tournament is constructed from the members of a suitable finite field by connecting pairs of elements that differ by a quadratic residue). In \cite{Gr}, Gr\"{u}nbaum proved that every tournament is anti-traceable unless it is isomorphic to one of $PT_{3},RT_{5}$ and $PT_{7}$. Rosenfeld in \cite{Ro} showed that for any tournament on $n\geq 9$ vertices, every vertex is the initial vertex of a Hamiltonian anti-directed path, strengthening the previous result. In general digraphs, Ghouila-Houri \cite{G} proved that every strong digraph $D$ on $n$ vertices with $\delta^{+}(D)+\delta^{-}(D)\ge n$ is Hamiltonian. Particularly, any digraph $D$ of $\delta^{0}(D)\geq n/2$ is Hamiltonian.  Woodall \cite{W} improved Ghouila-Houri's result by verifying that every strong digraph on $n\geq 2$ vertices with $d_D^+(u)+d^-_D(v)\ge n$ for every pair of nonadjacent vertices $u$ and $v$ in $D$ is Hamiltonian.

Recall that an oriented graph is a digraph without directed cycles of length two.
In oriented graphs, there exist a number of beautiful and significant results concerning Hamiltonian oriented cycles. A problem proposed by Thomassen \cite{T} is to determine $\delta^{0}(D)$ that guarantees an oriented graph to be Hamiltonian. H\"{a}ggkvist \cite{H} constructed a family of non-Hamiltonian oriented graphs $D$ with  $n$ vertices and $\delta^0(D)=(3n-5)/8$. Furthermore, H\"{a}ggkvist \cite{H} conjectured that any oriented graph $D$ with $n$ vertices and  $\delta^{0}(D)\geq(3n-4)/8$ is Hamiltonian. This conjecture was confirmed by Keevash, K\"{u}hn and Osthus \cite{KKO09} for sufficiently large $n$.
H\"{a}ggkvist \cite{H} also conjectured that any oriented graph with $n$ vertices and  $\delta(D)+\delta^+(D)+\delta^-(D)>(3n-3)/2$ is Hamiltonian. This conjecture was resolved by Kelly, K\"{u}hn and Osthus \cite{KKO08} approximately.
H\"{a}ggkvist and Thomason \cite{HT97} proposed a conjecture stating that for any positive $\epsilon$, the condition $\delta^{0}(D)\geq (3/8+\epsilon)n$ is sufficient to guarantee the existence of any orientation of a Hamiltonian cycle in any sufficiently large oriented graph $D$. The same authors subsequently demonstrated that the stronger condition $\delta^{0}(D)\geq (5/12+\epsilon)n$ is also adequate. In 2011, Kelly \cite{K} completely resolved this conjecture. Notably, this minimum semi-degree condition ensures that $D$ contains a cycle of every possible orientation and length.

Several interesting results have also been obtained for long oriented paths and long oriented cycles with specific lengths.
Jackson \cite{J} showed the following classical result.

\begin{theorem}[\cite{J}]\label{THM: jackson}
For any positive integer $k$, every oriented graph $D$ of $\delta^{0}(D)\geq k$ admits a directed path of length $2k$ when $|V(D)|>2k+2$. Furthermore, if $|V(D)|\le 2k+2$, then $D$ is Hamiltonian.
\end{theorem}

Zhang \cite{Z} refined this result by showing the following result.

\begin{theorem}[\cite{Z}]\label{THM:Zh}
Every oriented graph $D$ with $\delta^{0}(D)\geq k$ admits either a directed path of length at least $2k+2$ or a directed cycle of length at least $2k+1$.
\end{theorem}

Building upon these earlier results, Stein \cite{S} proposed the following conjecture.

\begin{conjecture}\textnormal{(\cite{S})} \label{conj}
For any positive integer $k$, every oriented graph $D$ of $\delta^{0}(D)>k/2$ contains every oriented path of length $k$.
\end{conjecture}

Conjecture~\ref{conj} is, in a sense, best possible due to the existence of disjoint union of $RT_{k}$.
Theorems~\ref{THM: jackson} and \ref{THM:Zh} showed that Conjecture~\ref{conj} holds for directed paths. The following is a list of known results concerning Conjecture~\ref{conj}.
\begin{itemize}
	\item[(1)] Kelly~\cite{K} proved that Conjecture~\ref{conj} is true for $k\geq 3n/4+o(n)$.
	
	\item[(2)] Jackson~\cite{J} proved that Conjecture~\ref{conj} holds for directed paths.
	
	\item[(3)]
	Stein and Z\'{a}rate-Guer\'{e}n~\cite{SZ} proved that for all $\eta\in (0, 1)$, there is $n_0$ such that for all $n\ge n_0$ and $k\ge \eta n$, every oriented graph
	$D$ on $n$ vertices with $\delta^0(D)>(1 +\eta)k/2$ contains every antidirected path with $k$ edges.
	
	\item[(4)] Klimo\v{s}ov\'{a} and Stein~\cite{KS} proved that for all $k\geq 3$, every oriented graph $D$ of $\delta^{0}(D)\geq (3k-2)/4$ contains each anti-directed path of length $k$.
	
	\item[(5)] Chen, Hou, and Zhou~\cite{CHZ} improved this lower bound to $(2k+1)/3$.
	
	\item[(6)] Skokan and Tyomkyn~\cite{ST} further improved the lower bound to $5k/8$.
\end{itemize}


For integers $s$ and $t$, define $P(s,t)$ to be an oriented path consisting of $s$ consecutive forward arcs followed by $t$ consecutive backward arcs.  Up to isomorphism, we have $P(s,t)=P(t,s)$.
Jackson's result (see (2)) can be rephrased as Conjecture~\ref{conj} holds true for $P(k,0)$. The main result of this paper is the following.


\begin{theorem} \label{THM: main}

For any integer $k$, every oriented graph $D$ with $\delta^{0}(D)>k/2$ contains a copy of $P(s,t)$ for all $s+t=k$.
\end{theorem}

Similarly, we define $Q(s,t)$ to be an oriented path of length $s+t$, consisting of $s$ backward arcs and followed by $t$ forward arcs. Also, $Q(s,t)$ is isomorphic to $Q(t,s)$. Notably, Theorem~\ref{THM: main} implies that every oriented graph $D$ with $\delta^{0}(D)>k/2$ also contains every $Q(s,t)$ for $s+t=k$. Here is the proof: consider the oriented graph  $\overleftarrow{D}$, obtained from $D$ by reversing the orientation of each arc. Then we have $\delta^{0}(\overleftarrow{D})>k/2$. By Theorem~\ref{THM: main}, $\overleftarrow{D}$ contains a $P(s,t)$, which is a $Q(s,t)$ included in $D$.

A \emph{block} of an oriented path $P$ is a maximal directed subpath in $P$. Consequently, an oriented path with two blocks is either isomorphic to a $P(s,t)$ or to a $Q(s,t)$ for some integers $s$ and $t$. As a corollary, we conclude that Conjecture~\ref{conj} is true for every oriented path with two blocks.

\begin{theorem} \label{THM: Main}

For any integer $k$, every oriented graph $D$ of $\delta^{0}(D)>k/2$ contains every oriented path of length $k$ with two blocks.

\end{theorem}

The rest of this paper is arranged as follows. In Section 2, we give the proof of Theorem~\ref{THM: main}. We conclude this paper with a remark in the last section.

\section{The proof of Theorem~\ref{THM: main}}

\noindent

We begin by introducing additional notation. Let $D=(V,A)$ be a digraph. The notation $u\rightarrow v$ stands for  $(u,v)\in A(D)$, while, $u\nrightarrow v$ indicates that $(u,v)\notin A(D)$.
The \emph{out-neighborhood} of a vertex $v$, denoted as $N^{+}_{D}(v)$, consists of all vertices $u\in V(D)$ such that $v\rightarrow u$, and the \emph{out-degree} of $v$, denoted as $d^{+}_{D}(v)$, is the cardinality of $N^{+}_{D}(v)$.
Analogously, the \emph{in-neighborhood} of a vertex $v$ is defined as $N^{-}_{D}(v)=\{u\in V(D)|\;u\rightarrow v\}$, and the \emph{in-degree} of $v$ is given by $d^{-}_{D}(v)=|N^{-}_{D}(v)|$. The vertices within $N^{+}_{D}(v)$ are referred to as \emph{out-neighbors}, while those in $N^{-}_{D}(v)$ are  considered \emph{in-neighbors} of $v$.
For integers $a$ and $b$, we use the notation $[a,b]=\{a, a+1, \dots, b\}$, and for convenience, we further denote $[k]=[1,k]$.


The following observation is straightforward yet highly valuable.

\begin{observation}\label{observation}

	Let $\ell$ be a positive integer and let $D$ be an oriented graph. If $P=v_{1}v_{2}\ldots v_{\ell}$ is a longest directed path in $D$, then 
	$N^{-}_{D}(v_{1}) \cup N^{+}_{D}(v_{\ell})\subseteq V(P)$.

\end{observation}

Now we are ready to give the proof of Theorem~\ref{THM: main}.

\vspace{0.3cm}

\noindent \emph{Proof of Theorem~\ref{THM: main}:}
Let $D$ be an oriented graph with $\delta^{0}(D)>k/2$. For $k\leq 4$, the statement can be directly verified. Hence, we assume $k\geq 5$ in the following. Suppose to the contrary that $D$ contains no path $P(s,t)$ for some $s$ and $t$ with $s+t=k$. By symmetry, we can assume that $s\geq t$. Therefore, $s\geq \lceil k/2\rceil$ and  $t\leq\lfloor k/2\rfloor$. Moreover, by Theorem~\ref{THM: jackson}, we may assume $t\geq 1$.

Let  $P=v_{1}v_{2}\ldots v_{\alpha}$ be a longest directed path in $D$.
By Observation~\ref{observation}, we have $N^{-}_{D}(v_{1}) \cup N^{+}_{D}(v_{\alpha})\subseteq V(P)$.
By Theorem~\ref{THM: jackson},   $\alpha\geq 2\delta^0(D)+1\ge k+2$.

Let $R=\{v_{s+1},v_{s+2},\ldots, v_{\alpha-t-1}\}$ and $B=\{v_{t+1},v_{t+2},\ldots,v_{\alpha-s-1}\}$. Observe that $|R|=|B|=\alpha-s-t-1\geq 1$ since $\alpha\geq k+2$ and $s+t=k$.

\begin{claim}\label{R}

\textnormal{(i)} No vertex in $R$ is adjacent to $v_{\alpha}$; \textnormal{(ii)} No vertex in $B$ is adjacent to $v_{\alpha}$.

\end{claim}
\begin{proof}

We only prove (i) as (ii) can be proved with similar discussion. Suppose the assertion is false. Then there must exist $i\in[s+1,\alpha-t-1]$ such that $v_{i}$ and $v_{\alpha}$ are adjacent.
If $v_{i}\rightarrow v_{\alpha}$, then$$v_{i-s+1}\rightarrow v_{i-s+2}\rightarrow \ldots \rightarrow v_{i}\rightarrow v_{\alpha}\leftarrow v_{\alpha-1}\leftarrow \ldots \leftarrow v_{\alpha-t}$$
is a path $P(s,t)$, a contradiction.
If $v_{\alpha}\rightarrow v_{i}$, then$$v_{i-s}\rightarrow v_{i-s+1}\rightarrow \ldots \rightarrow v_{i}\leftarrow v_{\alpha}\leftarrow v_{\alpha-1}\leftarrow \ldots \leftarrow v_{\alpha-t+1}$$
is a  $P(s,t)$, a contradiction too.
\end{proof}

From Claim~\ref{R}, we have the following technical statement.

\begin{claim}\label{CL: longestpathcontra}

	There exists no vertex $v$ in $D$ such that $v$ serves as the initial vertex of a longest directed path $P_1$ and also as the terminal vertex of another longest directed path $P_2$ with $V(P_1)=V(P_2)$.

\end{claim}
\begin{proof}

	Without loss of generality, assume $V(P_1)=V(P_2)=V(P)$.
	By Observation~\ref{observation}, we have $N^{+}_{D}(v)\cup N^{-}_{D}(v)\subseteq V(P)$.
	From Claim~\ref{R}, it is clear that $v$ is not adjacent to at least $\alpha-k-1$ vertices of $P$. Since $D$ is an oriented graph, it follows that $d^{+}_{D}(v)+ d^{-}_{D}(v)\leq k$. On the other hand, we observe that $d^{+}_{D}(v)+ d^{-}_{D}(v)\geq 2\delta^{0}(D)>k$. Consequently, we deduce that $k<d^{+}_{D}(v)+ d^{-}_{D}(v)\leq k$, which is a contradiction. 
\end{proof}

\begin{claim}\label{va}
$v_{\alpha}\nrightarrow v_{\alpha-t}$ and $v_{\alpha}\nrightarrow v_{\alpha-s}$.
\end{claim}
\begin{proof}
If $v_{\alpha}\rightarrow v_{\alpha-t}$, then$$v_{\alpha-k}\rightarrow v_{\alpha-k+1}\rightarrow \ldots \rightarrow v_{\alpha-t}\leftarrow v_{\alpha}\leftarrow v_{\alpha -1}\leftarrow \ldots \leftarrow v_{\alpha-t+1}$$
is a $P(s,t)$, a contradiction.
If $v_{\alpha}\rightarrow v_{\alpha-s}$, then$$v_{\alpha-s+1}\rightarrow v_{\alpha-s+2}\rightarrow \ldots \rightarrow v_{\alpha}\rightarrow v_{\alpha-s}\leftarrow v_{\alpha-s-1}\leftarrow \ldots \leftarrow v_{\alpha-k}$$
is a  $P(s,t)$, also a contradiction.
\end{proof}

Let $\gamma=\text{min}\{i\, |\,v_{\alpha}\rightarrow v_{i}\}$ and $\xi=\text{max}\{i\,|\, v_{i}\rightarrow v_{1}\}$.
Since $\delta^{0}(D)>k/2\geq t$ and $N^{+}_{D}(v_{\alpha})\subseteq V(P)$, by Claim~\ref{R}, it is easily seen that $\gamma\leq s$.

\begin{claim}\label{Rv}
$N^{-}_{D}(v_{1})\cap (R\cup\{v_{\alpha-t}\})=\emptyset$.
\end{claim}
\begin{proof}
Suppose not, there exists $v_{i}\in R\cup\{v_{\alpha-t}\}$ such that $v_{i}\rightarrow v_{1}$. Then
$$v_{i-s+\gamma}\rightarrow v_{i-s+\gamma+1}\rightarrow \ldots \rightarrow v_{i}\rightarrow v_{1}\rightarrow v_{2}\rightarrow \ldots \rightarrow v_{\gamma}\leftarrow v_{\alpha}\leftarrow v_{\alpha-1}\leftarrow \ldots \leftarrow v_{\alpha-t+1}$$
is a  $P(s,t)$, a contradiction.
\end{proof}

\begin{claim}\label{gammaxi}
$\gamma\geq 2$, i.e., $\xi\leq \alpha-1$.
\end{claim}
\begin{proof}
If not, then $\gamma=1$ and $P+(v_\alpha, v_1)$ is a directed cycle of length $\alpha$. By the maximality of $P$ and Observation~\ref{observation}, we have $N^{+}_{D}(v)\cup N^{-}_{D}(v)\subseteq V(P)$ for any $v\in V(P)$.
By Claim~\ref{R}, we have $N^{+}_{D}(v_{\alpha})\cup N^{-}_{D}(v_{\alpha})\subseteq V(P)\backslash (R\cup\{v_{\alpha}\})$. As $D$ is an oriented graph, $N^{+}_{D}(v_{\alpha})\cap N^{-}_{D}(v_{\alpha})=\emptyset$.
Therefore, we have  $$k<2\delta^{0}(D)\leq d^{+}_{D}(v_{\alpha})+ d^{-}_{D}(v_{\alpha})\leq \alpha-|R|-1=s+t=k,$$ this is impossible. The claim follows.
\end{proof}

For convenience, denote $U=\{v_{1},v_{2},\ldots,v_{\gamma}\}$ and $W=\{v_{\xi},v_{\xi+1},\ldots,v_{\alpha}\}$. In the subsequent proof, we divide the discussion into two cases: when $2\leq \gamma\leq t$ and when $\gamma>t$.

\vspace{5pt}
\noindent{\bf{Case A}. $2\leq \gamma\leq t$}

Recall that $N^{-}_{D}(v_{1})\subseteq V(P)$ and $d^{-}_{D}(v_{1})\geq \delta^{0}(D)>k/2\ge t$. It is obvious that $\xi\geq t+3$.

\begin{claim}\label{Bv}
$N^{-}_{D}(v_{1})\cap (B\cup\{v_{\alpha-s}\})=\emptyset$.
\end{claim}
\begin{proof}
On the contrary, suppose that there exists $v_{i}\in B\cup\{v_{\alpha-s}\}$ such that $v_{i}\rightarrow v_{1}$. It is straightforward to check that$$v_{\alpha-s+1}\rightarrow v_{\alpha-s+2}\rightarrow \ldots \rightarrow v_{\alpha}\rightarrow v_{\gamma}\leftarrow v_{\gamma-1}\leftarrow \ldots \leftarrow v_{1}\leftarrow v_{i}\leftarrow v_{i-1}\leftarrow \ldots \leftarrow v_{i-t+\gamma}$$
is a $P(s,t)$, a contradiction.
\end{proof}

By Claims~\ref{Bv} and~\ref{Rv}, we have $\xi\in [\alpha-s+1,s]\cup [\alpha-t+1,\alpha-1]$. Let $F=\{v_{\alpha-s},v_{\alpha-s+1},\ldots,v_{s}\}$. Set $X_{1}=\{v_{\gamma+1},v_{\gamma+2},\ldots,v_{t}\}$ and $Y_{1}=\{v_{\alpha-t},v_{\alpha-t+1},\ldots,v_{\xi-1}\}$.
Then $$V(P)=U\cup X_1\cup B\cup F\cup R\cup Y_1\cup W.$$
We divide the proof of Case A into two subcases: when $F=\emptyset$ and when $F\not=\emptyset$.

\vspace{0.2cm}
{\bf Case 1}. $F=\emptyset$
\vspace{0.2cm}

Then $\xi\in [\alpha-t+1,\alpha-1]$ and, clearly, $\alpha\geq 2s+1$. See Figure 1 for an illustration.

\begin{figure}[hbtp]
	\begin{center}	
		\begin{tikzpicture}[thick,scale=1, every node/.style={transform shape}]\label{fig_1}			
			\node (1) at (0,0) {$v_1$};
			\node (2) at (0.8,0) {$v_2$};
			\node (3) at (1.3,0) {$\ldots$};
			\node (4) at (2.25,0) {$v_\gamma$};
			\node (5) at (3.25,0) {$v_{\gamma+1}$};
			\node (6) at (3.9,0) {$\ldots$};
			\node (7) at (4.8,0) {$v_t$};
			\node (8) at (5.8,0) {$v_{t+1}$};
			\node (9) at (7,0) {$v_{t+2}$};
			\node (10) at (7.6,0) {$\ldots$};
			\node (11) at (8.5,0) {$v_s$};
			\node (12) at (9.5,0) {$v_{s+1}$};
			\node (13) at (10.1,0) {$\ldots$};
			\node (13) at (11.45,0) {$v_{\alpha-s-1}$};
			\node (14) at (12.95,0) {$v_{\alpha-s}$};
			\node (15) at (13.6,0) {$\ldots$};
			\node (16) at (14.9,0) {$v_{\alpha-t-1}$};
			\node (17) at (0.7,-1.5) {$v_{\alpha-t}$};
			\node (18) at (1.3,-1.5) {$\ldots$};
			\node (19) at (2.4,-1.5) {$v_{\xi-1}$};
			\node (20) at (3.4,-1.5) {$v_{\xi}$};
			\node (21) at (4.45,-1.5) {$v_{\xi+1}$};
			\node (22) at (5.05,-1.5) {$\ldots$};
			\node (23) at (6,-1.5) {$v_{\alpha}$};
			\draw[->, line width=0.5pt] (0.15,0) -- (0.6,0);
			\draw[->, line width=0.5pt] (1.55,0) -- (2,0);
			\draw[->, line width=0.5pt] (2.45,0) -- (2.85,0);
			\draw[->, line width=0.5pt] (4.15,0) -- (4.6,0);
			\draw[->, line width=0.5pt] (5.0,0) -- (5.4,0);
			\draw[->, line width=0.5pt] (6.15,0) -- (6.6,0);
			\draw[->, line width=0.5pt] (7.85,0) -- (8.3,0);
			\draw[->, line width=0.5pt] (8.7,0) -- (9.1,0);
			\draw[->, line width=0.5pt] (10.35,0) -- (10.8,0);
			\draw[->, line width=0.5pt] (12.05,0) -- (12.5,0);
			\draw[->, line width=0.5pt] (13.85,0) -- (14.3,0);
			\draw[->, line width=0.5pt] (-0.15,-1.5) -- (0.3,-1.5);
			\draw[->, line width=0.5pt] (1.55,-1.5) -- (2,-1.5);
			\draw[->, line width=0.5pt] (2.75,-1.5) -- (3.2,-1.5);
			\draw[->, line width=0.5pt] (3.6,-1.5) -- (4.05,-1.5);
			\draw[->, line width=0.5pt] (5.3,-1.5) -- (5.75,-1.5);
			\draw[green] (-0.25,-0.45) rectangle (2.5,0.5);
			\node[green] at (1.125,0.7) {$U$};
			\draw[green] (2.75,-0.45) rectangle (5,0.5);
			\node[green] at (3.775,0.7) {$X_1$};
			\draw[blue] (5.25,-0.65) rectangle (12.2,0.7);
			\node[blue] at (8.725,0.9) {$B$};
			\draw[red] (8.25,-0.35) rectangle (15.5,0.4);
			\node[red] at (13.725,0.6) {$R$};
			\draw[green] (0.2,-1.85) rectangle (2.8,-1.10);
			\node[green] at (1.5,-0.9) {$Y_1$};
			\draw[green] (3.05,-1.85) rectangle (6.25,-1.10);
			\node[green] at (4.65,-0.9) {$W$};
			\draw[line width=0.5pt, postaction={decorate}, decoration={markings, mark=at position 0.5 with {\arrow[line width=0.3pt, scale=1.7]{To}}}](6,-1.5) .. controls (5,0) and (3.25,-1) .. (2.25,0);
			\draw[line width=0.5pt, postaction={decorate}, decoration={markings, mark=at position 0.5 with {\arrow[line width=0.3pt, scale=1.7]{To}}}](3.4,-1.5) .. controls (3.4,0) and (0,-1) .. (0,0);
		\end{tikzpicture}

		\vspace{0.5em}
		{\small \textbf{Figure 1}}
	\end{center}
\end{figure}

\begin{claim}\label{desired01}
One of the following must occur:

\textnormal{(i)}. There are vertices $v_{\ell}\in X_{1}$ and $v_{\ell-2}\in X_{1}\cup \{v_{\gamma}\}$ such that $v_{\alpha}\rightarrow v_{\ell}$ and $v_{\ell-2}\rightarrow v_{1}$.

\textnormal{(ii)}. There are vertices $v_{m}\in Y_{1}$ and $v_{m-2}\in Y_{1}$ such that $v_{\alpha}\rightarrow v_{m}$ and $v_{m-2}\rightarrow v_{1}$.
\end{claim}
\begin{proof}
On the contrary, neither (i) nor (ii) occurs.
Let $X_{1}^{*}=N^{+}_{D}(v_{\alpha})\cap X_{1}$ and $Y_{1}^{*}=N^{+}_{D}(v_{\alpha})\cap Y_{1}$. Denote by $x_{1}^{*}=|X_{1}^{*}|$ and $y_{1}^{*}=|Y_{1}^{*}|$.
Since (i) does not occur, we have $|N^{-}_{D}(v_{1})\cap(X_{1}\cup \{v_{\gamma}\})|\le |X_1\cup\{v_\gamma\}|-(x_1^*-1)$.
By Claim~\ref{va}, $v_{\alpha}\nrightarrow v_{\alpha-t}$. Since (ii) does not occur, we have $|N^{-}_{D}(v_{1})\cap Y_{1}|\le|Y_1|-(y_1^*-1)$.
By  the maximality of $\xi$, Claims~\ref{Rv} and~\ref{Bv}, we obtain that
\begin{align}
|N^{-}_{D}(v_{1})|
&=|N^{-}_{D}(v_{1})\cap (U\backslash \{v_{\gamma}\})|+|N^{-}_{D}(v_{1})\cap(X_{1}\cup \{v_{\gamma}\})|+|N^{-}_{D}(v_{1})\cap (Y_{1}\cup \{v_{\xi}\})|\nonumber\\
&\leq \gamma-3+|X_{1}|+1-(x_{1}^{*}-1)+|Y_{1}|-(y_{1}^{*}-1)+1\nonumber\\
&=\gamma-3+t-\gamma+1-(x_{1}^{*}-1)+\xi-\alpha+t-(y_{1}^{*}-1)+1\nonumber\\
&=2t-\alpha-x_{1}^{*}-y_{1}^{*}+\xi+1,\label{EQ:e1}
\end{align}
and, by the minimality of $\gamma$ and Claim~\ref{R}, we have
\begin{align}
|N^{+}_{D}(v_{\alpha})|
&=|X_{1}^{*}\cup \{v_{\gamma}\}|+|Y_{1}^{*}|+|N^{+}_{D}(v_{\alpha})\cap W|\nonumber\\
&\leq x_{1}^{*}+1+y_{1}^{*}+\alpha-\xi-1\nonumber\\
&=\alpha+x_{1}^{*}+y_{1}^{*}-\xi.\label{EQ:e2}
\end{align}
Combining (\ref{EQ:e1}) and (\ref{EQ:e2}), we have
\begin{align*}
d^{-}_{D}(v_{1})+d^{+}_{D}(v_{\alpha})&\leq 2t+1\leq 2\lfloor k/2\rfloor+1.
\end{align*}
On the other hand, $d^{-}_{D}(v_{1})+d^{+}_{D}(v_{\alpha})\geq 2\delta^{0}(D)>k$. Therefore, we conclude that
$$k<d^{-}_{D}(v_{1})+d^{+}_{D}(v_{\alpha})\leq 2\lfloor k/2\rfloor+1.$$
If $k$ is odd, then we have $k<2\lfloor k/2\rfloor+1=k$, a contradiction. Now assume $k$ is even. Then we have $\delta^{0}(D)\ge k/2+1$. Thus we have $$k+2\leq 2\delta^{0}(D)\leq d^{-}_{D}(v_{1})+d^{+}_{D}(v_{\alpha})\leq 2\lfloor k/2\rfloor+1=k+1,$$
a contradiction too.
As a consequence, the claim follows.
\end{proof}



If (i) of Claim~\ref{desired01} occurs, i.e., there exist $v_{\ell}\in X_{1}$ and $v_{\ell-2}\in X_{1}\cup \{v_{\gamma}\}$ such that $v_{\alpha}\rightarrow v_{\ell}$ and $v_{\ell-2}\rightarrow v_{1}$, then both
$$v_{\ell-1}\rightarrow v_{\ell}\rightarrow \ldots \rightarrow v_{\alpha}\rightarrow v_{\gamma}\rightarrow v_{\gamma+1}\rightarrow \ldots \rightarrow v_{\ell-2}\rightarrow v_{1}\rightarrow v_{2}\rightarrow \ldots \rightarrow v_{\gamma-1}$$
and
$$v_{\xi+1}\rightarrow v_{\xi+2}\rightarrow \ldots \rightarrow v_{\alpha}\rightarrow v_{\ell}\rightarrow v_{\ell+1}\rightarrow \ldots \rightarrow v_{\xi}\rightarrow v_{1}\rightarrow v_{2}\rightarrow \ldots \rightarrow v_{\ell-1}$$
are longest directed paths in $D$.
However, this  contradicts Claim~\ref{CL: longestpathcontra}.

If (ii) of Claim~\ref{desired01} occurs, i.e., there exist $v_{m-2},v_{m}\in Y_{1}$ such that $v_{\alpha}\rightarrow v_{m}$ and $v_{m-2}\rightarrow v_{1}$, then both
$$v_{m-1}\rightarrow v_{m}\rightarrow \ldots \rightarrow v_{\alpha}\rightarrow v_{\gamma}\rightarrow v_{\gamma+1}\rightarrow \ldots \rightarrow v_{m-2}\rightarrow v_{1}\rightarrow v_{2}\rightarrow \ldots \rightarrow v_{\gamma-1}$$
and
$$v_{\xi+1}\rightarrow v_{\xi+2}\rightarrow \ldots \rightarrow v_{\alpha}\rightarrow v_{m}\rightarrow v_{m+1}\rightarrow \ldots \rightarrow v_{\xi}\rightarrow v_{1}\rightarrow v_{2}\rightarrow \ldots \rightarrow v_{m-1}$$
are longest directed paths in $D$, a contradiction  too.

\vspace{0.2cm}
{\bf Case 2}. $F\neq\emptyset$
\vspace{0.2cm}

It is obvious that $\alpha\leq 2s$. See Figure 2 for an illustration.

\begin{figure}[hbtp]
	\begin{center}	
		\begin{tikzpicture}[thick,scale=1, every node/.style={transform shape}]\label{fig_1}			
			\node (1) at (0,0) {$v_1$};
			\node (2) at (0.8,0) {$v_2$};
			\node (3) at (1.3,0) {$\ldots$};
			\node (4) at (2.25,0) {$v_\gamma$};
			\node (5) at (3.25,0) {$v_{\gamma+1}$};
			\node (6) at (3.9,0) {$\ldots$};
			\node (7) at (4.8,0) {$v_t$};
			\node (8) at (5.8,0) {$v_{t+1}$};
			\node (9) at (7,0) {$v_{t+2}$};
			\node (10) at (7.6,0) {$\ldots$};
			\node (11) at (8.9,0) {$v_{\alpha-s-1}$};
			\node (12) at (10.35,0) {$v_{\alpha-s}$};
			\node (13) at (11,0) {$\ldots$};
			\node (13) at (11.95,0) {$v_{s}$};
			\node (14) at (0.65,-1.5) {$v_{s+1}$};
			\node (15) at (1.9,-1.5) {$v_{s+2}$};
			\node (16) at (2.5,-1.5) {$\ldots$};
			\node (17) at (3.8,-1.5) {$v_{\alpha-t-1}$};
			\node (18) at (5.25,-1.5) {$v_{\alpha-t}$};
			\node (19) at (5.85,-1.5) {$\ldots$};
			\node (20) at (6.95,-1.5) {$v_{\xi-1}$};
			\node (21) at (7.95,-1.5) {$v_{\xi}$};
			\node (22) at (9,-1.5) {$v_{\xi+1}$};
			\node (23) at (9.6,-1.5) {$\ldots$};
			\node (24) at (10.6,-1.5) {$v_{\alpha}$};
			\draw[->, line width=0.5pt] (0.15,0) -- (0.6,0);
			\draw[->, line width=0.5pt] (1.55,0) -- (2,0);
			\draw[->, line width=0.5pt] (2.45,0) -- (2.85,0);
			\draw[->, line width=0.5pt] (4.15,0) -- (4.6,0);
			\draw[->, line width=0.5pt] (5.0,0) -- (5.4,0);
			\draw[->, line width=0.5pt] (6.15,0) -- (6.6,0);
			\draw[->, line width=0.5pt] (7.85,0) -- (8.3,0);
			\draw[->, line width=0.5pt] (9.45,0) -- (9.9,0);
			\draw[->, line width=0.5pt] (11.25,0) -- (11.7,0);
			\draw[->, line width=0.5pt] (4.4,-1.5) -- (4.85,-1.5);
			\draw[->, line width=0.5pt] (6.1,-1.5) -- (6.55,-1.5);
			\draw[->, line width=0.5pt] (7.3,-1.5) -- (7.75,-1.5);
			\draw[->, line width=0.5pt] (8.15,-1.5) -- (8.6,-1.5);
			\draw[->, line width=0.5pt] (9.85,-1.5) -- (10.3,-1.5);
			\draw[->, line width=0.5pt] (2.75,-1.5) -- (3.15,-1.5);
			\draw[->, line width=0.5pt] (1,-1.5) -- (1.45,-1.5);
			\draw[->, line width=0.5pt] (-0.2,-1.5) -- (0.25,-1.5);
			\draw[green] (-0.25,-0.45) rectangle (2.5,0.5);
			\node[green] at (1.125,0.7) {$U$};
			\draw[green] (2.75,-0.45) rectangle (5,0.5);
			\node[green] at (3.775,0.7) {$X_1$};
			\draw[blue] (5.25,-0.45) rectangle (9.55,0.5);
			\node[blue] at (7.4,0.7) {$B$};
			\draw[green] (9.75,-0.45) rectangle (12.25,0.5);
			\node[green] at (11,0.7) {$F$};
			\draw[red] (0.25,-1.85) rectangle (4.35,-1.10);
			\node[red] at (2.3,-0.9) {$R$};			
			\draw[green] (4.8,-1.85) rectangle (7.4,-1.10);
			\node[green] at (6.3,-0.9) {$Y_1$};
			\draw[green] (7.65,-1.85) rectangle (10.85,-1.10);
			\node[green] at (9.25,-0.9) {$W$};
			\draw[line width=0.5pt, postaction={decorate}, decoration={markings, mark=at position 0.5 with {\arrow[line width=0.3pt, scale=1.7]{To}}}](10.6,-1.5) .. controls (5,0) and (3.25,-1) .. (2.25,0);
			\draw[line width=0.5pt, postaction={decorate}, decoration={markings, mark=at position 0.5 with {\arrow[line width=0.3pt, scale=1.7]{To}}}](7.95,-1.5) .. controls (3.4,0) and (0,-1) .. (0,0);
		\end{tikzpicture}

		\vspace{0.5em}
		{\small \textbf{Figure 2}}
	\end{center}
\end{figure}

\begin{claim}\label{desired2}
One of the following must occur:

\textnormal{(i)}. There are vertices $v_{\ell}\in X_{1}$ and $v_{\ell-2}\in X_{1}\cup \{v_{\gamma}\}$ with $v_{\alpha}\rightarrow v_{\ell}$ and $v_{\ell-2}\rightarrow v_{1}$.

\textnormal{(ii)}. There are vertices $v_{m}, v_{m-2}\in F$ with $v_{\alpha}\rightarrow v_{m}$ and  $v_{m-2}\rightarrow v_{1}$.

\textnormal{(iii)}. There are vertices $v_{h}, v_{h-2}\in Y_{1}$ with $v_{\alpha}\rightarrow v_{h}$ and $v_{h-2}\rightarrow v_{1}$.

\end{claim}
\begin{proof}
Suppose that none of (i), (ii) and (iii) occurs. Let $X_{1}^{*}=N^{+}_{D}(v_{\alpha})\cap X_{1}$, $F^{*}=N^{+}_{D}(v_{\alpha})\cap F$ and $Y_{1}^{*}=N^{+}_{D}(v_{\alpha})\cap Y_{1}$. Denote by $x_{1}^{*}=|X_{1}^{*}|$, $f^{*}=|F^{*}|$ and $y_{1}^{*}=|Y_{1}^{*}|$. Since (i) does not occur, we have  $|N^-_D(v_{1})\cap (X_1\cup\{v_\gamma\})|\le |X_1\cup\{v_\gamma\}|-(x_1^*-1)$. By Claim~\ref{va}, we get $v_{\alpha}\nrightarrow v_{\alpha-t}$ and $v_{\alpha}\nrightarrow v_{\alpha-s}$. Since (ii) and (iii) do not occur, we have $|N^-_D(v_{1})\cap Y_1|\le |Y_1|-(y_{1}^{*}-1)$ and $|N^-_D(v_{1})\cap F|\le |Y_1|-(f^{*}-1)$.
Therefore, we have
\begin{align}
|N^{-}_{D}(v_{1})|
&=|N^{-}_{D}(v_{1})\cap (U\setminus\{v_\gamma\})|+|N^{-}_{D}(v_{1})\cap (X_{1}\cup\{v_\gamma\})|+|N^{-}_{D}(v_{1})\cap F|\nonumber\\
&\quad+|N^{-}_{D}(v_{1})\cap (Y_{1}\cup \{v_{\xi}\})|\nonumber\\
&\leq \gamma-3+|X_{1}|+1-(x_{1}^{*}-1)+|F|-(f^{*}-1)+|Y_{1}|-(y_{1}^{*}-1)+1\nonumber\\
&=\gamma-3+t-\gamma+1-(x_{1}^{*}-1)+2s-\alpha+1-(f^{*}-1)\nonumber\\
&\quad+\xi-\alpha+t-(y_{1}^{*}-1)+1\nonumber\\
&=2s+2t-2\alpha-x_{1}^{*}-f^{*}-y_{1}^{*}+\xi+3\nonumber\\
&=2k-2\alpha-x_{1}^{*}-f^{*}-y_{1}^{*}+\xi+3,\label{EQ:e4}
\end{align}
and
\begin{align}
|N^{+}_{D}(v_{\alpha})|
&=|X_{1}^{*}\cup \{v_{\gamma}\}|+|F^{*}|+|Y^{*}_{1}|+|N^{+}_{D}(v_{\alpha})\cap W|\nonumber\\
&\leq x_{1}^{*}+1+f^{*}+y_{1}^{*}+\alpha-\xi-1\nonumber\\
&=\alpha+x_{1}^{*}+f^{*}+y_{1}^{*}-\xi.\label{EQ:e5}
\end{align}
Combining $(\ref{EQ:e4})$ and $(\ref{EQ:e5})$, we have
\begin{align}
d^{-}_{D}(v_{1})+d^{+}_{D}(v_{\alpha})&\leq 2k-\alpha+3.\label{EQ:e6}
\end{align}
By $d^{-}_{D}(v_{1})+d^{+}_{D}(v_{\alpha})\geq 2\delta^{0}(D)>k$, we obtain that
\begin{equation}\label{EQ: e16}
	k<d^{-}_{D}(v_{1})+d^{+}_{D}(v_{\alpha})\leq 2k-\alpha+3.
\end{equation}
Recall that $\alpha \geq k+2$. Thus we have $\alpha=k+2$ by (\ref{EQ: e16}), and $d^{-}_{D}(v_{1})+d^{+}_{D}(v_{\alpha})=k+1$.
If $k$ is even, then we have $k+2\le2\delta^{0}(D)\le d^{-}_{D}(v_{1})+d^{+}_{D}(v_{\alpha})=k+1$, a contradiction.
Now assume that $k$ is odd. Then $\delta^{0}(D)\ge (k+1)/2$. Note that the longest directed path $P$ has length $\alpha-1=k+1\le 2\delta^{0}(D)$.
By Theorem~\ref{THM:Zh}, $D$ contains a directed cycle $C$ of length $2\delta^{0}(D)+1=k+2$. Choose $P$ as the longest directed path in $C$. This implies that $v_{\alpha}\rightarrow v_{1}$. However, this contradicts Claim~\ref{gammaxi}.

Therefore, the claim follows.
\end{proof}



If (i) of Claim~\ref{desired2} occurs, i.e., there exist $v_{\ell}\in X_{1}$ and $v_{\ell-2}\in X_{1}\cup \{v_{\gamma}\}$ such that $v_{\alpha}\rightarrow v_{\ell}$ and $v_{\ell-2}\rightarrow v_{1}$, then both
$$v_{\ell-1}\rightarrow v_{\ell}\rightarrow \ldots \rightarrow v_{\alpha}\rightarrow v_{\gamma}\rightarrow v_{\gamma+1}\rightarrow \ldots \rightarrow v_{\ell-2}\rightarrow v_{1}\rightarrow v_{2}\rightarrow \ldots \rightarrow v_{\gamma-1}$$
and
$$v_{\xi+1}\rightarrow v_{\xi+2}\rightarrow \ldots \rightarrow v_{\alpha}\rightarrow v_{\ell}\rightarrow v_{\ell+1}\rightarrow \ldots \rightarrow v_{\xi}\rightarrow v_{1}\rightarrow v_{2}\rightarrow \ldots \rightarrow v_{\ell-1}$$
are longest directed paths in $D$. We have the desired contradiction to Claim~\ref{CL: longestpathcontra}.


If (ii) of Claim~\ref{desired2} occurs, i.e., there exist $v_{m-2},v_{m}\in F$ such that $v_{\alpha}\rightarrow v_{m}$ and $v_{m-2}\rightarrow v_{1}$, then
 both
$$v_{m-1}\rightarrow v_{m}\rightarrow \ldots \rightarrow v_{\alpha}\rightarrow v_{\gamma}\rightarrow v_{\gamma+1}\rightarrow \ldots \rightarrow v_{m-2}\rightarrow v_{1}\rightarrow v_{2}\rightarrow \ldots \rightarrow v_{\gamma-1}$$
and
$$v_{\xi+1}\rightarrow v_{\xi+2}\rightarrow \ldots \rightarrow v_{\alpha}\rightarrow v_{m}\rightarrow v_{m+1}\rightarrow \ldots \rightarrow v_{\xi}\rightarrow v_{1}\rightarrow v_{2}\rightarrow \ldots \rightarrow v_{m-1}$$
are longest directed paths in $D$, a contradiction to Claim~\ref{CL: longestpathcontra} too.


If (iii) of Claim~\ref{desired2} occurs, i.e., there exist $v_{h-2},v_{h}\in Y_{1}$ such that $v_{\alpha}\rightarrow v_{h}$ and $v_{h-2}\rightarrow v_{1}$,
then
$$v_{h-1}\rightarrow v_{h}\rightarrow \ldots \rightarrow v_{\alpha}\rightarrow v_{\gamma}\rightarrow v_{\gamma+1}\rightarrow \ldots \rightarrow v_{h-2}\rightarrow v_{1}\rightarrow v_{2}\rightarrow \ldots \rightarrow v_{\gamma-1}$$
and
$$v_{\xi+1}\rightarrow v_{\xi+2}\rightarrow \ldots \rightarrow v_{\alpha}\rightarrow v_{h}\rightarrow v_{h+1}\rightarrow \ldots \rightarrow v_{\xi}\rightarrow v_{1}\rightarrow v_{2}\rightarrow \ldots \rightarrow v_{h-1}$$
are longest directed paths in $D$, the same contradiction as above.


\vspace{5pt}
\noindent{\bf{Case  B}. $\gamma >t$}

Then $t<\gamma\leq s$. Recall that $(N^{+}_{D}(v_{\alpha})\cup N^{-}_{D}(v_{\alpha}))\cap R=\emptyset$, $N^{-}_{D}(v_{1})\cap (R\cup \{v_{\alpha-t}\})=\emptyset$ and $v_{\alpha}\nrightarrow v_{\alpha-t}$.

\begin{claim}\label{xi}
$\xi\geq \gamma+4$.
\end{claim}
\begin{proof}
By the minimality of $\gamma$ and since $D$ is oriented, we have $d^{+}_{D}(v_{\alpha}) \subseteq \{v_{\gamma},v_{\gamma+1},\ldots,v_{\alpha-2}\}$.
Further, since $N^{+}_{D}(v_{\alpha})\cap (R\cup \{v_{\alpha-t}\})=\emptyset$, we deduce that
$$d^{+}_{D}(v_{\alpha})\leq \alpha-2-\gamma+1-(\alpha-s-t)=k-\gamma-1.$$
By $d^{-}_{D}(v_{1})\geq 2\delta^{0}(D)-d^{+}_{D}(v_{\alpha})>\gamma+1$, we have $\xi\geq \gamma+4$.
\end{proof}

Together with Claim~\ref{gammaxi} and $N^{-}_{D}(v_{1})\cap (R\cup \{v_{\alpha-t}\})=\emptyset$, we obtain $\xi\in [\gamma+4,s]\cup [\alpha-t+1,\alpha-1]$. Here to proceed with the proof, we will consider the following two cases.

\vspace{0.2cm}

{\bf Case 1}. $\gamma +4\leq \xi\leq s$
\vspace{0.2cm}

Let us denote $X=\{v_{\gamma+1},v_{\gamma+2},\ldots,v_{\xi-1}\}$ and $Y=\{v_{\alpha-t},v_{\alpha-t+1},\ldots,v_{\alpha}\}$. Write $X^{*}$ for $N^{+}_{D}(v_{\alpha})\cap X$ and let $x^{*}=|X^{*}|$. See Figure 3 for an illustration.

\begin{figure}[hbtp]
	\begin{center}	
		\begin{tikzpicture}[thick,scale=1, every node/.style={transform shape}]\label{fig_1}			
			\node (1) at (0,0) {$v_1$};
			\node (2) at (0.8,0) {$v_2$};
			\node (3) at (1.3,0) {$\ldots$};
			\node (4) at (2.25,0) {$v_\gamma$};
			\node (5) at (3.25,0) {$v_{\gamma+1}$};
			\node (6) at (3.9,0) {$\ldots$};
			\node (7) at (5,0) {$v_{\xi-1}$};
			\node (8) at (6,0) {$v_{\xi}$};
			\node (9) at (7,0) {$v_{\xi+1}$};
			\node (10) at (7.6,0) {$\ldots$};
			\node (11) at (8.5,0) {$v_{s}$};
			\node (12) at (9.55,0) {$v_{s+1}$};
			\node (13) at (10.15,0) {$\ldots$};
			\node (13) at (11.45,0) {$v_{\alpha-t-1}$};
			\node (14) at (12.9,0) {$v_{\alpha-t}$};
			\node (15) at (13.5,0) {$\ldots$};
			\node (17) at (14.45,0) {$v_{\alpha}$};
			\draw[->, line width=0.5pt] (0.15,0) -- (0.6,0);
			\draw[->, line width=0.5pt] (1.55,0) -- (2,0);
			\draw[->, line width=0.5pt] (2.45,0) -- (2.85,0);
			\draw[->, line width=0.5pt] (4.15,0) -- (4.6,0);
			\draw[->, line width=0.5pt] (5.35,0) -- (5.8,0);
			\draw[->, line width=0.5pt] (6.15,0) -- (6.6,0);
			\draw[->, line width=0.5pt] (7.85,0) -- (8.3,0);
			\draw[->, line width=0.5pt] (8.7,0) -- (9.15,0);
			\draw[->, line width=0.5pt] (10.4,0) -- (10.85,0);
			\draw[->, line width=0.5pt] (12,0) -- (12.45,0);
			\draw[->, line width=0.5pt] (13.75,0) -- (14.2,0);
			\draw[green] (-0.25,-0.45) rectangle (2.5,0.5);
			\node[green] at (1.125,0.7) {$U$};
			\draw[green] (2.75,-0.45) rectangle (5.4,0.5);
			\node[green] at (4.075,0.7) {$X$};
			\draw[red] (9,-0.45) rectangle (12.1,0.5);
			\node[red] at (10.55,0.7) {$R$};
			\draw[green] (12.25,-0.45) rectangle (14.75,0.5);
			\node[green] at (13.5,0.7) {$Y$};
			\draw[line width=0.5pt, postaction={decorate}, decoration={markings, mark=at position 0.5 with {\arrow[line width=0.3pt, scale=1.7]{To}}}](14.45,0) .. controls (8.35,1.5) .. (2.25,0);
			\draw[line width=0.5pt, postaction={decorate}, decoration={markings, mark=at position 0.5 with {\arrow[line width=0.3pt, scale=1.7]{To}}}](6,0) .. controls (3,-1) .. (0,0);
		\end{tikzpicture}
		
		\vspace{0.5em}
		{\small \textbf{Figure 3}}
	\end{center}
\end{figure}

\begin{claim}\label{desired3}
There are vertices $v_{\ell}\in X$ and $v_{\ell-2}\in X\cup \{v_{\gamma}\}$ with $v_{\alpha}\rightarrow v_{\ell}$ and $v_{\ell-2}\rightarrow v_{1}$.
\end{claim}
\begin{proof}
Suppose to the contrary that the claim is false. Then there are  at least $x^{*}-1$ vertices in $X\cup\{v_{\gamma}\}$ not belonging to $N_D^-(v_1)$. Thus we have
\begin{align}
|N^{-}_{D}(v_{1})|
&=|N^{-}_{D}(v_{1})\cap (U\backslash \{v_{\gamma}\})|+|N^{-}_{D}(v_{1})\cap (X\cup \{v_{\gamma},v_{\xi}\})|\nonumber\\
&\leq \gamma-3+|X|+1-(x^{*}-1)+1\nonumber\\
&=\gamma-3+\xi-\gamma-(x^{*}-1)+1\nonumber\\
&=\xi-x^{*}-1,\label{EQ:e7}
\end{align}
and
\begin{align}
|N^{+}_{D}(v_{\alpha})|
&=|X^{*}\cup \{v_{\gamma}\}|+|N^{+}_{D}(v_{\alpha})\cap Y|+|N^+_D(v_\alpha)\cap\{v_\xi,v_{\xi+1},\dots,v_s\}|\nonumber\\
&\leq x^{*}+1+t-2+s-\xi+1\nonumber\\
&=x^{*}+k-\xi.\label{EQ:e8}
\end{align}
By (\ref{EQ:e7}) and (\ref{EQ:e8}), we deduce that $d^{-}_{D}(v_{1})+d^{+}_{D}(v_{\alpha})\le k-1$.
Additionally, considering that $d^{-}_{D}(v_{1})+d^{+}_{D}(v_{\alpha})\geq 2\delta^{0}(D)>k$. We conclude that $k<d^{-}_{D}(v_{1})+d^{+}_{D}(v_{\alpha})\leq k-1$, a contradiction.
\end{proof}

By Claim~\ref{desired3}, there exist $v_{\ell}\in X$ and $v_{\ell-2}\in X\cup \{v_{\gamma}\}$ such that $v_{\alpha}\rightarrow v_{\ell}$ and $v_{\ell-2}\rightarrow v_{1}$, then both
$$v_{\ell-1}\rightarrow v_{\ell}\rightarrow \ldots \rightarrow v_{\alpha}\rightarrow v_{\gamma}\rightarrow v_{\gamma+1}\rightarrow \ldots \rightarrow v_{\ell-2}\rightarrow v_{1}\rightarrow v_{2}\rightarrow \ldots \rightarrow v_{\gamma-1}$$
and
$$v_{\xi+1}\rightarrow v_{\xi+2}\rightarrow \ldots \rightarrow v_{\alpha}\rightarrow v_{\ell}\rightarrow v_{\ell+1}\rightarrow \ldots \rightarrow v_{\xi}\rightarrow v_{1}\rightarrow v_{2}\rightarrow \ldots \rightarrow v_{\ell-1}$$
are longest directed paths in $D$ satisfying that $v_{\ell-1}$ is the initial vertex of one path and the terminal vertex of the other. This contradicts Claim~\ref{CL: longestpathcontra}.


\vspace{0.2cm}

{\bf Case 2}. $\alpha-t+1\leq \xi\leq \alpha-1$.
\vspace{0.2cm}

Let $X_{2}=\{v_{\gamma+1},v_{\gamma+2},\ldots,v_{s}\}$ and $Y_{2}=\{v_{\alpha-t},v_{\alpha-t+1},\ldots,v_{\xi-1}\}$.
Then $V(P)$ has the partition $U\cup X_2\cup R\cup Y_2\cup W$. See Figure 4 for an illustration.

\begin{figure}[hbtp]
	\begin{center}	
		\begin{tikzpicture}[thick,scale=1, every node/.style={transform shape}]\label{fig_1}			
			\node (1) at (0,0) {$v_1$};
			\node (2) at (0.8,0) {$v_2$};
			\node (3) at (1.3,0) {$\ldots$};
			\node (4) at (2.25,0) {$v_\gamma$};
			\node (5) at (3.25,0) {$v_{\gamma+1}$};
			\node (6) at (3.9,0) {$\ldots$};
			\node (7) at (4.8,0) {$v_s$};
			\node (8) at (5.8,0) {$v_{s+1}$};
			\node (9) at (7,0) {$v_{s+2}$};
			\node (10) at (7.6,0) {$\ldots$};
			\node (11) at (8.9,0) {$v_{\alpha-t-1}$};
			\node (12) at (10.35,0) {$v_{\alpha-t}$};			
			\node (13) at (11.8,0) {$v_{\alpha-t+1}$};
			\node (13) at (13.65,0) {$v_{\xi-1}$};
			\node (16) at (12.55,0) {$\ldots$};
			\node (21) at (0.5,-1.5) {$v_{\xi}$};
			\node (22) at (1.55,-1.5) {$v_{\xi+1}$};
			\node (23) at (2.15,-1.5) {$\ldots$};
			\node (24) at (3.15,-1.5) {$v_{\alpha}$};
			\draw[->, line width=0.5pt] (0.15,0) -- (0.6,0);
			\draw[->, line width=0.5pt] (1.55,0) -- (2,0);
			\draw[->, line width=0.5pt] (2.45,0) -- (2.85,0);
			\draw[->, line width=0.5pt] (4.15,0) -- (4.6,0);
			\draw[->, line width=0.5pt] (5.0,0) -- (5.4,0);
			\draw[->, line width=0.5pt] (6.15,0) -- (6.6,0);
			\draw[->, line width=0.5pt] (7.85,0) -- (8.3,0);
			\draw[->, line width=0.5pt] (9.45,0) -- (9.9,0);
			\draw[->, line width=0.5pt] (10.75,0) -- (11.2,0);
			\draw[->, line width=0.5pt] (12.8,0) -- (13.25,0);
			\draw[->, line width=0.5pt] (-0.15,-1.5) -- (0.3,-1.5);
			\draw[->, line width=0.5pt] (0.7,-1.5) -- (1.15,-1.5);
			\draw[->, line width=0.5pt] (2.4,-1.5) -- (2.85,-1.5);
			\draw[green] (-0.25,-0.45) rectangle (2.5,0.5);
			\node[green] at (1.125,0.7) {$U$};
			\draw[green] (2.75,-0.45) rectangle (5,0.5);
			\node[green] at (3.775,0.7) {$X_2$};
			\draw[red] (5.25,-0.45) rectangle (9.55,0.5);
			\node[red] at (7.4,0.7) {$R$};
			\draw[green] (9.75,-0.45) rectangle (14.15,0.5);
			\node[green] at (11.95,0.7) {$Y_2$};
			\draw[green] (0.2,-1.85) rectangle (3.4,-1.10);
			\node[green] at (1.8,-0.9) {$W$};
			\draw[line width=0.5pt, postaction={decorate}, decoration={markings, mark=at position 0.5 with {\arrow[line width=0.3pt, scale=1.7]{To}}}](3.15,-1.5) -- (2.25,0);
			\draw[line width=0.5pt, postaction={decorate}, decoration={markings, mark=at position 0.5 with {\arrow[line width=0.3pt, scale=1.7]{To}}}](0.5,-1.5) -- (0,0);
		\end{tikzpicture}
		
		\vspace{0.5em}
		{\small \textbf{Figure 4}}
	\end{center}
\end{figure}

\begin{claim}\label{desired1}
One of the following must occur:

\textnormal{(i)}. There are vertices $v_{\ell}\in X_{2}$ and $v_{\ell-2}\in X_{2}\cup \{v_{\gamma}\}$ with $v_{\alpha}\rightarrow v_{\ell}$ and $v_{\ell-2}\rightarrow v_{1}$.

\textnormal{(ii)}. There are vertices $v_{m}\in Y_{2}$ and $v_{m-2}\in Y_{2}$ with $v_{\alpha}\rightarrow v_{m}$ and $v_{m-2}\rightarrow v_{1}$.
\end{claim}
\begin{proof}
Suppose, on the contrary, that neither (i) nor (ii) occurs. Let $X_{2}^{*}=N^{+}_{D}(v_{\alpha})\cap X_{2}$ and $Y_{2}^{*}=N^{+}_{D}(v_{\alpha})\cap Y_{2}$. Denote by $x_{2}^{*}=|X_{2}^{*}|$ and $y_{2}^{*}=|Y_{2}^{*}|$. Then there are  at least $x_{2}^{*}-1$ vertices in $X_{2}\cup\{v_{\gamma}\}$ not belonging to $N_D^-(v_1)$, and there are at least $y_{2}^{*}-1$ vertices in $Y_{2}$ not belonging to $N_D^-(v_1)$ because $v_{\alpha}\nrightarrow v_{\alpha-t}$.
Consequently, we obtain that
\begin{align}
|N^{-}_{D}(v_{1})|
&=|N^{-}_{D}(v_{1})\cap (U\backslash \{v_{\gamma}\})|+|N^{-}_{D}(v_{1})\cap(X_{2}\cup \{v_{\gamma}\})|+|N^{-}_{D}(v_{1})\cap (Y_{2}\cup \{v_{\xi}\})|\nonumber\\
&\leq \gamma-3+|X_{2}|+1-(x_{2}^{*}-1)+|Y_{2}|-(y_{2}^{*}-1)+1\nonumber\\
&=\gamma-3+s-\gamma+1-(x_{2}^{*}-1)+\xi-\alpha+t-(y_{2}^{*}-1)+1\nonumber\\
&=k-\alpha-x_{1}^{*}-y_{1}^{*}+\xi+1,\label{EQ:e10}
\end{align}
and
\begin{align}
|N^{+}_{D}(v_{\alpha})|
&=|X_{2}^{*}\cup \{v_{\gamma}\}|+|Y_{2}^{*}|+|N^{+}_{D}(v_{\alpha})\cap W|\nonumber\\
&\leq x_{2}^{*}+1+y_{2}^{*}+\alpha-\xi-1\nonumber\\
&=\alpha+x_{2}^{*}+y_{2}^{*}-\xi.\label{EQ:e11}
\end{align}
Combining  (\ref{EQ:e10}) and (\ref{EQ:e11}), we have
\begin{align*}
d^{-}_{D}(v_{1})+d^{+}_{D}(v_{\alpha})&\leq k+1.\label{EQ:e12}
\end{align*}
However,  $$k<2\delta^{0}(D)\leq d^{-}_{D}(v_{1})+d^{+}_{D}(v_{\alpha})\leq k+1$$
forces that $d^{-}_{D}(v_{1})+d^{+}_{D}(v_{\alpha})=k+1$ and thus $d^{-}_{D}(v_{1})=k-\alpha-x_{1}^{*}-y_{1}^{*}+\xi+1$ and $d^{+}_{D}(v_{\alpha})=\alpha+x_{2}^{*}+y_{2}^{*}-\xi$. Particularly, we have $|N^{-}_{D}(v_{1})\cap Y_{2}|=|Y_{2}|-(y_{2}^{*}-1)$ and $|N^{+}_{D}(v_{\alpha})\cap W|=\alpha-\xi-1$. $|N^{-}_{D}(v_{1})\cap Y_{2}|=|Y_{2}|-(y_{2}^{*}-1)$ implies that $v_{\xi-2}\rightarrow v_{1}$, and  $|N^{+}_{D}(v_{\alpha})\cap W|=\alpha-\xi-1$ implies that  $v_{\alpha}\rightarrow v_{\xi}$.
Then
$$v_{\xi-1}\rightarrow v_{\xi}\rightarrow \ldots \rightarrow v_{\alpha}\rightarrow v_{\gamma}\rightarrow v_{\gamma+1}\rightarrow \ldots \rightarrow v_{\xi-2}\rightarrow v_{1}\rightarrow v_{2}\rightarrow \ldots \rightarrow v_{\gamma-1}$$
and
$$v_{\xi+1}\rightarrow v_{\xi+2}\rightarrow \ldots \rightarrow v_{\alpha}\rightarrow v_{\xi}\rightarrow v_{1}\rightarrow v_{2}\rightarrow \ldots \rightarrow v_{\xi-1}$$
are longest directed paths in $D$ satisfying that $v_{\xi-1}$ is the initial vertex of one path and the terminal vertex of the other. This contradicts Claim~\ref{CL: longestpathcontra}.


Therefore, the claim follows.
\end{proof}

With the same discussion as the above, if either (i) or (ii) of Claim~\ref{desired1} occurs, we have a vertex $v\in V(P)$ such that $v$ is the initial vertex of a longest directed path $P_1$ and the terminal vertex of another longest directed path $P_2$ with $V(P_1)=V(P_2)=V(P)$. This contradicts  Claim~\ref{CL: longestpathcontra}.


\vspace{0.2cm}

This completes our proof of the theorem.\hfill$\square$

\section{Remarks}

\noindent

Jackson~\cite{J} has demonstrated that Conjecture~\ref{conj} holds true for every oriented path with one block. In this paper, we show that Conjecture~\ref{conj} remains valid for every oriented path with two blocks. The natural question arises: Can we extend the confirmation of Conjecture~\ref{conj} to oriented paths with three or more blocks using similar techniques.

\subsection*{Acknowledgements}

This work was supported by the National Key Research and Development Program of China (2023YFA1010200), the National Natural Science Foundation of China (Nos. 12071453, 12471336), and the Innovation Program for Quantum Science and Technology (2021ZD0302902).

\vskip 3mm
\end{document}